\documentclass[12pt]{article}
\usepackage{epsfig,latexsym,amsfonts,amsmath,amssymb}

\usepackage{color}

\usepackage[notref,notcite]{showkeys}  % TO BE ERASED IN THE FINAL FORM

\def\textsf{\bf}  %%%%% here!!!!

\newtheorem{Theorem}{Theorem}
\newtheorem{Lemma}%[Theorem]
{Lemma}
\newtheorem{Proposition}%[Theorem]
{Proposition}
\newtheorem{Corollary}%[Theorem]
{Corollary}
\newtheorem{Remark}%[Theorem]
{Remark}

\def\proof{\noindent{\textbf{Proof. }}}
\def\R{\mathbb R}

\def\SoOm{{H^s(\Omega)}}
\def\SoOmNull{H^s_0(\Omega)}
\def\tildeSobolev{\widetilde{H}^s(\Omega)}
\def\SoOmsi{{H^\sigma(\Omega)}}
\def\SoOmNullsi{H^\sigma_0(\Omega)}
\def\tildeSobolevsi{\widetilde{H}^\sigma(\Omega)}
\def\SoOmmin{{H^{-\sigma}(\Omega)}}
\def\tildeSobolevmin{\widetilde{H}^{-\sigma}(\Omega)}
\def\SoRn{H^s(\mathbb R^n)}

\def\div{{\rm div}}

\begin{document}

\title %{\vspace{-10mm}
{On fractional Laplacians -- 2}

\author{Roberta Musina\footnote{Dipartimento di Matematica ed Informatica, Universit\`a di Udine,
via delle Scienze, 206 -- 33100 Udine, Italy. Email: {roberta.musina@uniud.it}. 
{Partially supported by Miur-PRIN 2009WRJ3W7-001 ``Fenomeni di concentrazione e {pro\-ble\-mi} di analisi geometrica''.}}~ and
Alexander I. Nazarov\footnote{St.Petersburg Department of Steklov Institute, Fontanka, 27, St.Petersburg, 191023, Russia
and St.Petersburg State University, 
Universitetskii pr. 28, St.Petersburg, 198504, Russia. E-mail: {al.il.nazarov@gmail.com}. Supported by RFBR grant 14-01-00534 and by
St.Petersburg University grant 6.38.670.2013.}
}

\date{}

\maketitle

\footnotesize

\noindent
{\bf Abstract.} The present paper is the natural evolution of arXiv:1308.3606. For $s>-1$ we compare two natural types of fractional Laplacians 
$(-\Delta)^s$, namely, the ``Navier'' and the ``Dirichlet'' ones. As a main tool, we give the ``dual'' Caffarelli--Silvestre and Stinga--Torrea
characterizations of these operators for $s\in(-1,0)$.

%\normalsize

%\bigskip\bigskip

\section{Introduction}

Recall that the Sobolev space $\SoRn=W^s_2(\mathbb{R}^n)$, $s\in\mathbb R$, is the space of distributions $u\in{\cal S}'(\mathbb R^n)$
with finite norm
$$
\|u\|_s^2=\int\limits_{\mathbb R^n}\left(1+|\xi|^2\right)^s|\mathcal Fu(\xi)|^2~\!d\xi,
$$
see for instance Section~2.3.3 of the monograph \cite{Tr}. Here ${\cal F}$ denotes
the Fourier transform
$$
\mathcal F{u}(\xi)=\frac{1}{(2\pi)^{n/2}}\int\limits_{\mathbb R^n} e^{-i~\!\!\xi\cdot x}u(x)~\!dx.
$$

For arbitrary $s\in\mathbb R$ we define fractional Laplacian in $\mathbb R^n$ by the quadratic 
form 
$$
Q_s[u]=((-\Delta)^su,u):=\int\limits_{\mathbb R^n}|\xi|^{2s}|{\cal F}u(\xi)|^2 d\xi,
$$
with domain
$$
{\rm Dom}(Q_s)=%\begin{cases}
%\{u\in L_2(\mathbb{R}^n)\,:\, Q_s[u]<\infty\}, & {\rm if} \ s\ge0;\\
\{u\in {\cal S}'(\mathbb{R}^n) \,:\, Q_s[u]<\infty\}.%, & {\rm if} \ s<0.
%\end{cases}
$$

Let $\Omega$ be a bounded and smooth domain in $\R^n$. We put
$$
\SoOm=\left\{u\big|_{\Omega}\,:\,u\in\SoRn\right\},
$$
see \cite[Sec.~4.2.1]{Tr}  and the extension theorem in \cite[Sec.~4.2.3]{Tr}.

Also we introduce the space
$$\tildeSobolev=\{u\in \SoRn\,:\,{\rm supp}\, u\subset\overline{\Omega}\}.
$$
By Theorem 4.3.2/1 \cite{Tr}, for $s-\frac{1}{2}\notin\mathbb{Z}$ this space coincides with 
$\SoOmNull$ that is the closure of ${\cal C}^{\infty}_0(\Omega)$ in $\SoOm$ while for 
$s-\frac{1}{2}\in\mathbb{Z}$ one has $\tildeSobolev\subsetneq\SoOmNull$. Moreover, 
${\cal C}^{\infty}_0(\Omega)$ is dense in $u\in \tildeSobolev$. \medskip

We introduce the ``Dirichlet'' fractional Laplacian in $\Omega$ (denoted by $(-\Delta_{\Omega})^s_D$) as the restriction of $(-\Delta)^s$. 
%$$
%Q_{s,\Omega}^D[u]=((-\Delta_{\Omega})^s_Du,u):=\int\limits_{\mathbb R^n}|\xi|^{2s}|{\cal F}u(\xi)|^2 d\xi~,
%$$
The domain of its quadratic form is
$$
{\rm Dom}(Q_{s,\Omega}^D)=\{u\in {\rm Dom}(Q_s) \,:\, {\rm supp}\, u\subset\overline{\Omega}\}.
$$

Also we define the ``Navier'' fractional Laplacian as $s$-th power of the conventional
Dirichlet Laplacian in the sense of spectral theory. Its quadratic form reads
$$
Q_{s,\Omega}^N[u]=((-\Delta_{\Omega})^s_Nu,u):=\sum\nolimits_j\lambda_j^s\cdot|(u,\varphi_j)|^2.
$$
Here, $\lambda_j$ and $\varphi_j$ are eigenvalues and eigenfunctions of the Dirichlet Laplacian in $\Omega$, respectively,
and ${\rm Dom}(Q_{s,\Omega}^N)$ consists of distributions in $\Omega$ such that $Q_{s,\Omega}^N[u]<\infty$.

It is well known that for $s=1$ these operators coincide: $(-\Delta_{\Omega})_N=(-\Delta_{\Omega})_D$.
We emphasize that, in contrast to $(-\Delta_{\Omega})^s_N$, the operator $(-\Delta_{\Omega})^s_D$ is not the $s$-th power of the 
Dirichlet Laplacian for $s\ne1$. In particular, $(-\Delta_{\Omega})^{-s}_D$ is not inverse to $(-\Delta_{\Omega})^s_D$.

%\blue{
%In addition, \red{(this is the former Lemma 1)} it holds that 
%$\tildeSobolev\subseteq\widetilde H^s_N(\Omega)$ and
%\begin{equation}
%\label{eq:spaces}
%\tildeSobolev=\widetilde H^s_N(\Omega)\quad\text{if and only if~ $s<\displaystyle\frac{3}{2}$}.
%\end{equation}
%This fact is well known for integer orders s; the general case  follows immediately
%from \cite[Theorem 1.17.1/1]{Tr} and \cite[Theorem 4.3.2/1]{Tr}. }

\medskip

The present paper is the natural evolution of \cite{FL}, where we compared the operators $(-\Delta_{\Omega})^s_D$ and $(-\Delta_{\Omega})^s_N$
for $0<s<1$. In the first result we extend Theorem 2 of \cite{FL}.

\begin{Theorem}\label{T:main1}
Let $s>-1$, $s\notin\mathbb N_0$. Then for $u\in {\rm Dom}(Q_{s,\Omega}^D)$, $u\not\equiv0$, the following relations hold:
\begin{eqnarray}
Q_{s,\Omega}^N[u] > Q_{s,\Omega}^D[u], & {\rm if} & 2k<s<2k+1,\ \ k\in\mathbb N_0;
\label{corol1}\\
Q_{s,\Omega}^N[u] < Q_{s,\Omega}^D[u], & {\rm if} & 2k-1<s<2k,\ \ k\in\mathbb N_0.
\label{corol2}
\end{eqnarray}
\end{Theorem}

Next, we take into account the role of dilations in $\R^n$. We denote by $F(\Omega)$
the class of smooth and bounded domains containing $\Omega$. If $\Omega'\in F(\Omega)$,
then any $u\in {\rm Dom}(Q_{s,\Omega}^D)$ can be regarded as a function in
${\rm Dom}(Q_{s,\Omega'}^D)$, and the corresponding form $Q_{s,\Omega'}^D[u]$ does not change.
In contrast, the form $Q_{s,\Omega'}^N[u]$ does
depend on $\Omega'\supset\Omega$. However, roughly speaking, the difference between these quadratic
forms disappears as $\Omega'\to \R^n$.

\begin{Theorem}
\label{T:new}
Let $s>-1$. Then for $u\in {\rm Dom}(Q_{s,\Omega}^D)$ the following facts hold:
\begin{eqnarray}
Q_{s,\Omega}^D[u]=\inf_{\Omega'\in F(\Omega)}Q_{s,\Omega'}^N[u], & {\rm if} & 2k<s<2k+1,\ \ k\in\mathbb N_0;
\label{new1}\\
Q_{s,\Omega}^D[u]=\sup_{\Omega'\in F(\Omega)}Q_{s,\Omega'}^N[u], & {\rm if} & 2k-1<s<2k,\ \ k\in\mathbb N_0.
\label{new2}
\end{eqnarray}
\end{Theorem}

For $-1<s<0$ we also obtain a pointwise comparison result reverse to the case $0<s<1$ (compare with \cite[Theorem 1]{FL}).

\begin{Theorem} 
\label{T:main2}
Let $-1<s<0$, and let $f\in {\rm Dom}(Q_{s,\Omega}^D)$, $f\ge0$ in the sense of distributions, $f\not\equiv0$. Then 
the following relation holds:
\begin{equation}
(-\Delta_{\Omega})^s_Nf < (-\Delta_{\Omega})^s_Df.
\label{-pos_pres}
\end{equation}
\end{Theorem}

Actually, fractional Laplacians of orders $s\in (-1,0)$ play a crucial role in our arguments. 
In Section \ref{negat} we give a variational characterization of %negative Laplacians
these operators, ``dual'' to variational characterization of fractional Laplacians of orders 
$s\in (0,1)$ given in \cite{CfSi} and \cite{ST}. Theorems \ref{T:main1}--\ref{T:main2} are proved 
in Section \ref{compar}. 

Note that our statements hold in more general setting. Let $\Omega$ be a bounded and smooth domain 
in a complete smooth Riemannian manifold ${\cal M}$. Denote by $(-\Delta_{\Omega})^s_N$ and 
$(-\Delta_{\Omega})^s_D$, respectively, the $s$-th power of the Dirichlet Laplacian in $\Omega$ 
and the restriction of $s$-th power of the Dirichlet Laplacian in ${\cal M}$ to the set 
of functions supported in $\Omega$. Then proofs of Theorems \ref{T:main1}--\ref{T:main2} (and of Theorem 1
in \cite{FL} as well) run with minimal changes.

% \begin{Remark}
%Theorems \ref{T:main1} and \ref{T:main2} hold indeed in  more general settings.
%
%Let $\mathcal M$ be a complete Riemannian manifold with empty or smooth boundary, 
%and let $\Omega$ be a smooth domain in $\mathcal M$. We define standard Dirichlet Laplacians $-\Delta_M$
%on $M$ and $-\Delta_\Omega$ on $\Omega$. For $s>-1$ we put
%$$
%(-\Delta_\Omega)^s_D:=(-\Delta_M)^s\big|_{\Omega}~,\quad (-\Delta_\Omega)^s_N:=(-\Delta_\Omega)^s~\!,
%$$
%where the right-hand sides have to be intended in the sense of spectral theory. 
%The proofs run with without essential changes, thanks to the general extension results by
%Stinga and Torrea \cite{ST}. 
%
%\medskip
%
%Theorem \ref{T:new} can be extended as well. Notice however that 
%if $\mathcal M$ is compact without boundary, then it can happen that
%\begin{equation}
%\label{eq:arc}
%Q_{s,\Omega}^D[u]\neq \lim_{\Omega'\to \mathcal M\atop \Omega'\in F(\Omega)}Q_{s,\Omega'}^N[u].
%\end{equation}
%For instance,
%take $\mathcal M=\mathbb S^1$ to be the unit circle and fix a nontrivial $u\in C^\infty_0(\Omega)$.
%For any $z_0\in \mathbb S^1\setminus\Omega$
%one can find a sequence  of  arcs $\Omega'\supset\Omega$ such that $\Omega'\to \mathbb S^1\setminus\{z_0\}$
%and $Q_{s,\Omega'}^N[u]\to Q_{s,\mathbb S^1\setminus\{z_0\}}^N[u]$. Therefore, if $s>(1/2,1)$ then
%the  limit in the right hand side of (\ref{eq:arc}) does not exist, and
%$$
%Q_{s,\Omega}^D[u]< \liminf_{\Omega'\to \mathcal M\atop \Omega'\in F(\Omega)}Q_{s,\Omega'}^N[u].
%$$
%
%\end{Remark}

\section{Fractional Laplacians of negative orders}\label{negat}

First, we recall some facts from the classical monograph \cite{Tr} about the spaces $\SoOm$ and $\tildeSobolev$.

\begin{Proposition}\label{Hs} (a particular case of \cite[Theorem 4.3.2/1]{Tr}).

\begin{enumerate}
 \item If $0< \sigma<\frac 12$ then $\tildeSobolevsi=\SoOmNullsi=\SoOmsi$;
 \item If $\sigma=\frac 12$ then $\tildeSobolevsi$ is dense in $\SoOmsi=\SoOmNullsi$;
 \item If $\frac 12<\sigma< 1$ then $\tildeSobolevsi=\SoOmNullsi$ is a subspace of $\SoOmsi$.
 \end{enumerate}
\end{Proposition}

\begin{Proposition}\label{duality} (a particular case of \cite[Theorem 2.10.5/1]{Tr}).

For any $\sigma\in\mathbb R$ $(\tildeSobolevsi)'=\SoOmmin$.
\end{Proposition}
As an immediate consequence we obtain

\begin{Corollary}\label{H-s}
 
\begin{enumerate}
 \item If $0< \sigma<\frac 12$ then $\tildeSobolevmin=\SoOmmin$;
 \item If $\sigma=\frac 12$ then $\tildeSobolevmin$ is dense in $\SoOmmin$;
 \item If $\frac 12<\sigma< 1$ then $\SoOmmin$ is a subspace of $\tildeSobolevmin$.
 \end{enumerate}
\end{Corollary}

\begin{Remark}\label{1D}
In the one-dimensional case, for $\frac 12<\sigma< 1$ the codimension of $\SoOmmin$ in $\tildeSobolevmin$ equals $2$ since
the same is codimension of $\tildeSobolevsi$ in $\SoOmsi$. 
\end{Remark}

The next statement gives explicit description of domains of quadratic forms under consideration.

\begin{Lemma}\label{domain}
Let $0<\sigma <1$. Then
\begin{enumerate}
 \item  ${\rm Dom}(Q_{-\sigma ,\Omega}^N)=\SoOmmin$;
 \item  ${\rm Dom}(Q_{-\sigma ,\Omega}^D)=\tildeSobolevmin$ if $n\ge2$ or $\sigma <\frac 12$;
 \item  ${\rm Dom}(Q_{\sigma ,\Omega}^D)=\{u\in\tildeSobolevmin\,:\,{\cal F}u(0)=0\}$ if $n=1$ and $\sigma \ge\frac 12$.
\end{enumerate}
\end{Lemma}

\proof
The first statement follows from the relation ${\rm Dom}(Q_{\sigma ,\Omega}^N)=\tildeSobolevsi$, see, e.g.,
\cite[Theorems 1.15.3 and 4.3.2/2]{Tr}, and from Proposition \ref{duality}.

The second and the third statements follow directly from definition of $\tildeSobolevmin$, if we take into account
that ${\cal F}u$ is a smooth function.
\hfill$\square$\medskip

By Lemma \ref{domain} and Corollary \ref{H-s}, for $0<\sigma \le\frac 12$ we have ${\rm Dom}(Q_{-\sigma ,\Omega}^D)\subseteq{\rm Dom}(Q_{-\sigma ,\Omega}^N)$
(even ${\rm Dom}(Q_{-\sigma ,\Omega}^D)={\rm Dom}(Q_{-\sigma ,\Omega}^N)$ if $0<\sigma <\frac 12$). In the case $\frac 12<\sigma < 1$, ${\rm Dom}(Q_{-\sigma ,\Omega}^N)$
is a subspace of ${\rm Dom}(Q_{-\sigma ,\Omega}^D)$ (for $n=1$ this follows from Remark \ref{1D}). However, for arbitrary 
$f\in{\rm Dom}(Q_{-\sigma ,\Omega}^D)$ we can consider $f$ as a functional on $\SoOmsi$, put $\widetilde f=f|_{\tildeSobolevsi}\in{\rm Dom}(Q_{-\sigma ,\Omega}^N)$
and define $Q_{-\sigma ,\Omega}^N[f]:=Q_{-\sigma ,\Omega}^N[\widetilde f]$. \medskip

Next, we recall that in the paper \cite{CfSi} the fractional Laplacian of order
$\sigma \in(0,1)$ in  $\mathbb R^n$ was connected with the so-called {\it harmonic extension in $n+2-2\sigma $ dimensions} 
and with generalized Dirichlet-to-Neumann map (see also \cite{CT} for the case $\sigma =\frac 12$). 
In particular, for any $u\in\tildeSobolevsi$ the function $w_\sigma ^D(x,y)$ minimizing the weighted Dirichlet integral
$$
{\cal E}_\sigma ^D(w)=\int\limits_0^\infty\!\int\limits_{\mathbb{R}^n} y^{1-2\sigma }|\nabla w(x,y)|^2\,dxdy
$$
over the set
$${\cal W}_\sigma ^D(u)=\Big\{w(x,y)\,:\,{\cal E}_\sigma ^D(w)<\infty~,\ \ w\big|_{y=0}=u\Big\},
$$
satisfies
\begin{equation}
Q_{\sigma ,\Omega}^D[u]=\frac {C_\sigma }{2\sigma }\cdot {\cal E}_\sigma ^D(w_\sigma ^D),
\label{quad_D}
\end{equation}
where the constant $C_\sigma $ is given by
$$
C_\sigma :=\frac{4^\sigma \Gamma(1+\sigma )}{\Gamma(1-\sigma )}.
$$
Moreover, $w_\sigma ^D(x,y)$ is the solution of the BVP
$$-\div (y^{1-2\sigma }\nabla w)=0\quad \mbox{in}\quad \mathbb R^n\times\mathbb R_+;\qquad w\big|_{y=0}=u,
$$
and for sufficiently smooth $u$
\begin{equation}
(-\Delta)^\sigma u(x)=-\frac {C_\sigma }{2\sigma }\cdot\lim\limits_{y\to0^+} 
y^{1-2\sigma }\partial_yw_\sigma ^D(x,y),\qquad x\in\mathbb R^n
\label{extension_D}
\end{equation}
(we recall that $(-\Delta_{\Omega})^\sigma _Du=(-\Delta)^\sigma  u\big|_{\Omega}$).

In \cite{ST} this approach was developed in quite general situation. 
In particular, it was shown that for any $u\in\tildeSobolevsi$
the function $w_\sigma ^N(x,y)$ minimizing the energy integral
\begin{equation*}
{\cal E}_\sigma ^N(w)=\int\limits_0^\infty\!\int\limits_{\Omega} y^{1-2\sigma }|\nabla w(x,y)|^2\,dxdy
\label{energy_N}
\end{equation*}
over the set 
$$
{\cal W}_{\sigma ,\Omega}^N(u)=\{w(x,y)\in{\cal W}_\sigma ^D(u)\,:\,w\big|_{x\in\partial\Omega}=0\},
$$
satisfies
%and (\ref{extension_N}) implies
\begin{equation}
Q_{\sigma ,\Omega}^N[u]=\frac {C_\sigma }{2\sigma }\cdot {\cal E}_\sigma ^N(w_\sigma ^N).
\label{quad_N}
\end{equation}
Moreover, $w_\sigma ^N(x,y)$ is the solution of the BVP
\begin{equation}
-\div (y^{1-2\sigma }\nabla w)=0\quad \mbox{in}\quad \Omega\times\mathbb R_+;\qquad w\big|_{y=0}=u;
\qquad w\big|_{x\in\partial\Omega}=0,
\label{eq:ST}
\end{equation}
and for sufficiently smooth $u$ it turns out that
\begin{equation}
(-\Delta_{\Omega})^\sigma _Nu(x)=-\frac {C_\sigma }{2\sigma }\cdot\lim\limits_{y\to0^+} 
y^{1-2\sigma }\partial_yw_\sigma ^N(x,y).
\label{extension_N}
\end{equation}

In a similar way, negative fractional Laplacians are connected with generalized Neu\-mann-to-Dirichlet map. 
Namely, let $u\in {\rm Dom}(Q_{-\sigma ,\Omega}^D)$. We consider 
the problem of minimizing the functional 
$$
\widetilde{\cal E}_{-\sigma }^D(w)=%\int\limits_0^\infty\!\int\limits_{\mathbb{R}^n} y^{1-2\sigma }|\nabla w(x,y)|^2\,dxdy
{\cal E}_\sigma ^D(w)\,-\,2\,\big\langle u,w\big|_{y=0}\big\rangle
$$
over the set ${\cal W}_{-\sigma }^D$, that is closure of smooth functions on $\mathbb R^n\times\bar{\mathbb R}_+$ with bounded
support, with respect to ${\cal E}_\sigma ^D(\cdot)$. We recall that by Lemma 1 $u$ can be considered as a compactly supported
functional on $H^\sigma(\R^n)$, and thus the duality $\big\langle u,w\big|_{y=0}\big\rangle$ is well defined by the result of \cite{CfSi}.

First, let $n>2\sigma $ (this is a restriction only for $n=1$). We claim that the Hardy type inequality
\begin{equation}
{\cal E}_\sigma ^D(w)\ge \Big(\frac {n-2\sigma}{2}\Big)^2 \int\limits_0^\infty\!\int\limits_{\mathbb{R}^n} y^{1-2\sigma }\,\frac {w^2(x,y)}{r^2}\,dxdy
\label{eq:H}
\end{equation}
holds for $w\in{\cal W}_{-\sigma }^D$ (here $r^2=|x|^2+y^2$). Indeed, for a smooth function $w$ with bounded
support we consider the restriction of $w$ to arbitrary ray in $\mathbb R^n\times\mathbb R_+$ and write down the classical Hardy inequality
$$
\int\limits_0^\infty r^{n+1-2\sigma }w_r^2\,dr\ge \Big(\frac {n-2\sigma}{2}\Big)^2  \int\limits_0^\infty r^{n-1-2\sigma }w^2\,dr.
$$
We multiply it by $\big(\frac yr\big)^{1-2\sigma}$, integrate over unit hemisphere in $\mathbb R^{n+1}$, and the claim follows.

By (\ref{eq:H}), a non-zero constant cannot be approximated by compactly supported functions. 
Thus, the minimizer of $\widetilde{\cal E}_{-\sigma }^D$ is determined uniquely. Denote it by $w_{-\sigma }^D(x,y)$.
Then formulae (\ref{quad_D}) and (\ref{extension_D}) imply relations
\begin{equation}
Q_{-\sigma,\Omega }^D[u]=-\frac {2\sigma }{C_\sigma }\cdot \widetilde{\cal E}_{-\sigma }^D(w_{-\sigma }^D);\quad
(-\Delta_{\Omega})^{-\sigma }_Du(x)=\frac {2\sigma }{C_\sigma }\,w_{-\sigma }^D(x,0),\ \ x\in\Omega,
\label{-D}
\end{equation}
that give the ``dual'' Caffarelli--Silvestre characterization of $(-\Delta_{\Omega})^{-\sigma }_D$.

In case $n=1\le 2\sigma $ the above argument needs some modification. Namely, the minimizer $w_{-\sigma }^D(x,y)$ in this case
is defined up to an additive constant. However, by Lemma \ref{domain} we have
$${\cal F}u(0)\equiv\big\langle u,{\bf 1}\big\rangle=0.
$$
Therefore, $\widetilde{\cal E}_{-\sigma }^D(w_{-\sigma }^D)$ does not depend on the choice of the constant, and the first relation in (\ref{-D}) holds.
The second equality in (\ref{-D}) also holds if we choose the constant such that $w_{-\sigma }^D(x,0)\to0$ as $|x|\to\infty$.%\medskip

\begin{Remark}
Note that for sufficiently smooth $u$ the function $w_{-\sigma }^D$ solves the Neumann problem
\begin{equation}
-\div (y^{1-2\sigma }\nabla w)=0\quad \mbox{in}\quad \mathbb R^n\times\mathbb R_+;\qquad \lim\limits_{y\to0^+} 
y^{1-2\sigma }\partial_yw=-u.
\label{eq:-CS}
\end{equation}
\end{Remark}

Analogously, formulae (\ref{quad_N}) and (\ref{extension_N}) imply the ``dual'' Stinga--Torrea characterization of  $(-\Delta_{\Omega})^{-\sigma }_N$.
Namely, the function $w_{-\sigma }^N(x,y)$ minimizing the functional
\begin{equation*}
\widetilde{\cal E}_{-\sigma }^N(w)=%\int\limits_0^\infty\!\int\limits_{\Omega} y^{1-2\sigma }|\nabla w(x,y)|^2\,dxdy
{\cal E}_\sigma ^N(w)\,-\,2\,\big\langle u,w\big|_{y=0}\big\rangle
\end{equation*}
over the set 
$$
{\cal W}_{-\sigma ,\Omega}^N(u)=\{w(x,y)\in{\cal W}_{-\sigma }^D\,:\,w\big|_{x\notin\Omega}=0\},
$$
satisfies
\begin{equation}
Q_{-\sigma,\Omega }^N[u]=-\frac {2\sigma }{C_\sigma }\cdot \widetilde{\cal E}_{-\sigma }^N(w_{-\sigma }^N);\qquad
(-\Delta_{\Omega})^{-\sigma }_Nu(x)=\frac {2\sigma }{C_\sigma }\,w_{-\sigma }^N(x,0).
\label{-N}
\end{equation}

\begin{Remark}\label{ext}
Formula (\ref{-N}) shows that $w_{-\sigma }^N$ is the Stinga--Torrea extension of 
$\frac {C_\sigma }{2\sigma }(-\Delta_{\Omega})^{-\sigma }_Nu$. Similarly, from (\ref{-D}) 
we conclude that $w_{-\sigma }^D$ is the  Caffarelli--Silvestre extension of 
$\frac {C_\sigma }{2\sigma }(-\Delta)^{-\sigma }u$ but not of 
$\frac {C_\sigma }{2\sigma }(-\Delta_{\Omega})^{-\sigma }_Du$. This is due to the fact, already 
noticed in the introduction, that $(-\Delta_{\Omega})^{-\sigma }_D$ is not the inverse of 
$(-\Delta_{\Omega})^{\sigma }_D$.
\end{Remark}

\begin{Remark}
The representation of $(-\Delta)^{-\sigma}u$ via solution of the problem (\ref{eq:-CS}) was used in
\cite{CbS}. Similar representation of $(-\Delta_{\Omega})^{-\sigma}_Nu$ via solution of corresponding
mixed boundary value problem was used earlier in \cite{CDDS}. However, variational characterizations
of negative fractional Laplacians (the first parts of formulae (\ref{-D}) and (\ref{-N})) which play 
key role in what follows, are given for the first time.
\end{Remark}

\section{Proofs of main theorems}\label{compar}

We start by recalling an auxiliary result.

\begin{Lemma}\label{domain1}
Let $s>1$. Then 
$$\gathered
{\rm Dom}(Q_{s,\Omega}^D)=\tildeSobolev={\rm Dom}(Q_{s,\Omega}^N)\quad{\rm for}\quad s<3/2;\\
{\rm Dom}(Q_{s,\Omega}^D)=\tildeSobolev\subsetneq{\rm Dom}(Q_{s,\Omega}^N)\quad{\rm for}\quad s\ge3/2.
\endgathered
$$
\end{Lemma}

\proof 
For $Q_{s,\Omega}^D$ the conclusion follows directly from its definition. For $Q_{s,\Omega}^N$
this fact is well known for $s\in\mathbb N$; in general case it follows immediately
from \cite[Theorem 1.17.1/1]{Tr} and \cite[Theorem 4.3.2/1]{Tr}.\hfill$\square$\medskip

\noindent{\textbf{Proof of Theorem \ref{T:main1}.}} We split the proof in three parts.

{\bf 1}. Let $0<s<1$. Then the relation (\ref{corol1}) is proved in \cite[Theorem 2]{FL}.\medskip
 
{\bf 2}. Let $-1<s<0$. We define $\sigma=-s\in(0,1)$ and construct extensions $w_{-\sigma }^D$ and $w_{-\sigma }^N$ as described in Section \ref{negat}.
 
%If we assume a function $w\in{\cal W}_{-\sigma ,\Omega}^N$ extended by zero to $(\mathbb{R}^n\setminus\Omega)\times\mathbb{R}_+$ 
We evidently have ${\cal W}_{-\sigma ,\Omega}^N\subset{\cal W}_{-\sigma }^D$ and 
$\widetilde{\cal E}_{-\sigma }^N=\widetilde{\cal E}_{-\sigma }^D\big|_{{\cal W}_{-\sigma ,\Omega}^N}$. Therefore, (\ref{-D}) and (\ref{-N}) provide
$$Q_{s,\Omega}^N[u]=-\frac {2\sigma }{C_\sigma }\cdot \inf\limits_{w\in{\cal W}_{-\sigma ,\Omega}^N}\widetilde{\cal E}_{-\sigma }^N(w)
\le -\frac {2\sigma }{C_\sigma }\cdot \inf\limits_{w\in{\cal W}_{-\sigma }^D}\widetilde{\cal E}_{-\sigma }^D(w)=Q_{s,\Omega}^D[u].
$$
To complete the proof, we observe that for $u\not\equiv0$ the function $w_{-\sigma }^N$ cannot be a solution
of the homogeneous equation in (\ref{eq:ST}) in the whole half-space, %$\mathbb R^n\times\mathbb R_+$ 
since such a solution is analytic in the half-space. Thus, it cannot provide 
$\inf\limits_{w\in{\cal W}_{-\sigma }^D}\widetilde{\cal E}_{-\sigma }^D(w)$, and (\ref{corol2}) follows.\medskip

{\bf 3}. Now let $s>1$, $s\notin\mathbb N$. We put $k=\lfloor\frac {s-1}2\rfloor$ and define for $u\in\tildeSobolev$
$$
v=(-\Delta)^ku\in \widetilde H^{s-2k}(\Omega),\qquad s-2k\in (-1,0)\cup(0,1).
$$ Note that $v\not\equiv0$ if $u\not\equiv0$. Then we have
$$
%\aligned
Q_{s,\Omega}^N[u]=Q_{s-2k,\Omega}^N[v],\qquad Q_{s,\Omega}^D[u]=Q_{s-2k,\Omega}^D[u],
%\endaligned
$$
and the conclusion follows from cases 1 and 2.
\hfill$\square$\medskip

\noindent{\textbf{Proof of Theorem \ref{T:new}.}} Here we again distinguish three cases.

{\bf 1}. Let $0<s<1$. Then the relation (\ref{new1}) is proved in \cite[Theorem 3]{FL}.\medskip

{\bf 2}. Let $-1<s<0$. We define $\sigma=-s\in(0,1)$ and proceed similarly to the proof of \cite[Theorem 3]{FL}. 
It is sufficient to prove the statement for $u\in {\cal C}^\infty_0(\Omega)$.

For $\Omega'\supset\Omega$ we have ${\cal W}_{-\sigma ,\Omega}^N\subset {\cal W}_{-\sigma ,\Omega'}^N$. 
%(we recall that $w\in{\cal W}_{-\sigma ,\Omega}^N$ are assumed extended by zero to $(\mathbb{R}^n\setminus\Omega)\times\mathbb{R}_+$).
By (\ref{-N}), the quadratic form $Q_{s,\Omega}^N[u]$ is monotone increasing with respect to the domain inclusion.
Taking (\ref{corol2}) into account, we obtain 
\begin{equation}
\label{eq:monotone}
Q_{s,\Omega}^D[u]> Q_{s,\Omega'}^N[u]\ge Q_{s,\Omega}^N[u].
\end{equation}

Denote by $w=w^D_{-\sigma }$ the Caffarelli--Silvestre extension of $\frac {C_\sigma }{2\sigma }(-\Delta)^{-\sigma }u$,  
described in Section \ref{negat}.
 %Formula (\ref{-D}) implies that the quantity
%$$
%\int\limits_0^\infty \frac{1}{r}\Big\{r\int\limits_0^\infty\int\limits_{{\mathbb S}_{r}} y^{1-2\sigma }|\partial_y w(x,y)|^2
%~\!d{\mathbb S}_r(x)dy\Big\} 
%dr=\int\limits_0^\infty\irn y^{1-2\sigma }|\partial_y w(x,y)|^2~\!dxdy
%%<\infty,
%$$
%is finite (here $\mathbb S_r$ is the sphere of radius $r$ in $\mathbb R^n$).
%Since the function $r\mapsto r^{-1}$
%is not integrable at $\infty$, there exists a sequence $r_h\to\infty$ such that
%the balls $B_{r_h}$ contain $\Omega$ and
%\begin{equation}
%\label{eq:rh2}
%r_h\int\limits_0^\infty\int\limits_{{\mathbb S}_{r_h}} y^{1-2\sigma }|\partial_y w(x,y)|^2
%~\!d{\mathbb S}_{r_h}(x)dy \to 0.
%\end{equation}
Next, for any $y\ge 0$ let $\phi_R(\cdot,y)$ be the harmonic extension of
$w(\cdot,y)$ on the ball $B_R$, that is,
$$
%\begin{cases}
-\Delta \phi_R(\cdot,y)=0\quad\text{in} \ B_R;\qquad%\\
\phi_h(\cdot,y)=w(\cdot,y)\quad\text{on} \ \partial B_R~\!.
%\end{cases}
$$
Finally, for $x\in B_R$ and $y\ge0$ we put
$$
w_R(x,y)=
w(x,y)-\phi_R(x,y)~\!.
$$

It is shown in the proof of \cite[Theorem 3]{FL} that there exists a sequence $R_h\to\infty$ such that
$$
{\cal E}^N_\sigma (w_{R_h})\le {\cal E}^D_\sigma (w)+o(1).
$$
Further, since $(-\Delta_{\Omega})^{-\sigma }_Du$ vanishes at infinity, for any multi-index $\beta$ we evidently have $D^\beta\phi_{R_h}(\cdot,0)\to 0$
locally uniformly as $R_h\to\infty$. This gives $\big\langle u,\phi_{R_h}(\cdot,0)\big\rangle=o(1)$, and we obtain 
by (\ref{-D}) and (\ref{-N})
\begin{equation}
\label{eq:tesi}
\aligned
Q_{s,B_{R_h}}^N[u]
&\ge -\frac {2\sigma }{C_\sigma }\cdot \widetilde{\cal E}^N_{-\sigma }(w_{R_h})\\
&\ge -\frac {2\sigma }{C_\sigma }\cdot\widetilde{\cal E}^D_\sigma (w)-o(1)= Q_{s,\Omega}^D[u]-o(1)~\!.
\endaligned
\end{equation}
The relation (\ref{new1}) readily follows by comparing (\ref{eq:monotone}) and (\ref{eq:tesi}).\medskip

{\bf 3}. For $s>1$, $s\notin\mathbb N$, the conclusion follows from cases 1 and 2 just as in the proof of Theorem \ref{T:main1}.
\hfill$\square$\medskip

\begin{Remark}
\label{R:new}
Assume that $0\in\Omega$ and put $\alpha\Omega=\{\alpha x\,:\,x\in\Omega\}$.
Thanks to (\ref{eq:monotone}), the proof above shows indeed that
$$
Q^D_{s,\Omega}[u]=\lim_{\alpha\to\infty}Q^N_{s,\alpha\Omega}[u]
\quad\text{for any}\quad u\in\tildeSobolev.
$$
Now put $u_\alpha(x)=\alpha^{\frac{n-2s}{2}}u(\alpha x)$. Then the scaling shows that
$$
Q^D_{s,\Omega}[u_\alpha]\equiv Q^D_{s,\Omega}[u]=\lim_{\alpha\to\infty}Q^N_{s,\Omega}[u_\alpha]
\quad\text{for any}\quad u\in\tildeSobolev.
$$
\end{Remark}

\noindent{\textbf{Proof of Theorem \ref{T:main2}.}}
First, let $f\in{\cal C}^{\infty}_0(\Omega)$. We define $\sigma=-s\in(0,1)$ and construct extensions $w_{-\sigma }^D$ and $w_{-\sigma }^N$ 
described in Section \ref{negat}. 
%Note that they solve the corresponding BVPs for the Euler--Lagrange equations to functionals
%$\widetilde{\cal E}_{-\sigma }^D$ and $\widetilde{\cal E}_{-\sigma }^N$, respectively. 
Making the change of the variable $t=y^{2\sigma }$, we rewrite the BVP (\ref{eq:-CS}) for $w_{-\sigma }^D(x,t)$ as follows:
\begin{equation}
\Delta_xw_{-\sigma }^D+4\sigma ^2t^{\frac{2\sigma -1}\sigma }\partial^2_{tt}w_{-\sigma }^D=0\quad\mbox{in}\quad \mathbb R^n\times\mathbb R_+;
\qquad \partial_tw_{-\sigma }^D\big|_{t=0}=-\frac {f}{2\sigma }.
\label{BVP}
\end{equation}
Since $w_{-\sigma }^D$ vanishes at infinity, $w_{-\sigma }^D(x,t)>0$ for $t>0$ by the maximum principle. Moreover,
by \cite[Theorem 1.4]{ABMMZ} (the boundary point lemma) we have $w_{-\sigma }^D(x,0)>0$.

Further, the function $w_{-\sigma }^N$ satisfies the equalities (\ref{BVP}) in $\Omega\times\mathbb R_+$. Since
$w_{-\sigma }^N\big|_{x\notin\Omega}=0$, we infer that the function 
$$
W(x,t):=w_{-\sigma }^D(x,t)-w_{-\sigma }^N(x,t)
$$ 
meets the following relations:
\begin{equation*}
\Delta_xW+4\sigma ^2t^{\frac{2\sigma -1}\sigma }\partial^2_{tt}W=0\quad\mbox{in}\quad \Omega\times\mathbb R_+;
\qquad \partial_tW\big|_{t=0}=0;\qquad W\big|_{x\in\partial\Omega}>0.
\label{Maz}
\end{equation*}
Again, \cite[Theorem 1.4]{ABMMZ} gives $W(x,0)>0$, which gives (\ref{-pos_pres}) in view of (\ref{-D}) and (\ref{-N}).

For $f\in \tildeSobolev$ the statement holds by approximation argument.\hfill$\square$\medskip

%\section{A more general setting}\label{gener}

\end{document}